\documentclass[10pt,reqno]{amsart}
\usepackage{amsfonts,amsmath}
\usepackage{amssymb}
\usepackage{amsthm}
\usepackage{eucal}
 \newcommand{\R}{\mathbb R}

 \newcommand{\C}{\mathbb C}

 \newcommand{\p}{\partial}

 \renewcommand{\a}{a}

 \newcommand{\Chi}{{\bf \chi}}

\renewcommand{\b}{b}
\newcommand{\secao}[1]{\section{#1}\setcounter{equation}{0}}

\newtheorem{theorem}{Theorem}[section]
\newtheorem{proposition}[theorem]{Proposition}
 \newtheorem{remark}[theorem]{Remark}
\newtheorem{lemma}[theorem]{Lemma}

 \def\beqs{\begin{eqnarray*}}
 \def\enqs{\end{eqnarray*}}
 \def\beq{\begin{eqnarray}}
 \def\enq{\end{eqnarray}}
\begin{document}
\title[Schr\"odinger equation]{Sharp global well-posedness for a higher
order Schr\"odinger equation}
\setlength{\baselineskip}{1.3\baselineskip}
\author[X. Carvajal]{Xavier Carvajal}
\address{IMECC-UNICAMP, Campinas, 13081-970, Brazil}
\email{carvajal@ime.unicamp.br}

\thanks{The author was supported by FAPESP under grant No. 2004/07189-2.}
%%%%%%%%%%%%%%%%%%
\begin{abstract} Using the theory of almost conserved energies and the ``I-method'' developed by Colliander, Keel, Staffilani, Takaoka and Tao, we prove that the initial value problem for a
 higher order Schr\"odinger equation is globally well-posed in Sobolev spaces of order $s>1/4$. This result is sharp.
\end{abstract}
 \maketitle
\secao{Introduction}
 In this paper we will describe a sharp result of global
 well-posedness for solutions of the initial value problem (IVP)
\begin{equation}\label{1.1}
\begin{cases}
\p_t u+ia\,\p^2_xu+b\,\p^3_x u+ic\,|u|^{2}u+d\,
|u|^2 \p_x u+e\, u^2\p_x\bar{u}=  0, \quad x,t \in \R,\\
u(x,0)  =  \varphi(x),
\end{cases}
\end{equation}
where $u$ is a complex valued function and $a, b, c, d$ and $e$ are real 
parameters  with $be\neq 0$. %and $c= a(d-e)/(3b)$.

This model was proposed by Hasegawa and Kodama in \cite{H-K, Ko}
to describe the nonlinear propagation of pulses in optical fibers. In
literature, this model is called as a higher order nonlinear
Schr\"odinger equation or also Airy-Schr\"odinger equation.

We consider the following gauge transformation
\begin{equation}\label{GT}
v(x,t)=\exp\Big(i\lambda\,x+i( a\,\lambda^2-2b\,\lambda^3) t\Big)\,
u(x+(2a\, \lambda -3b \lambda^2)t,t),
\end{equation}
then, $u$ solves (\ref{1.1}) if and only if $v$
satisfies the IVP
\begin{equation}\label{1.9}
\begin{cases}
\p_t v+i(a-3\lambda\, b)\p_x^2 v+b\,\p_x^3 v+i(c-\lambda(d-e))\,|v|^2 v+d\,|v|^2\p_x v+ e\,v^2\p_x \bar v=0,\\
v(x,0)  = \exp(i\lambda\,x)\,
u(x,0).
\end{cases}
\end{equation}
Thus, if we take $\lambda=a/3b$ in (\ref{GT}) and $c=(d-e)a/3b$, then the function
\begin{equation}\label{1.20}
v(x,t)=\exp\Big(i\frac{a}{3b}\,x+i\frac{a^3}{27b^2}\,t\Big)\,
u(x+\frac{a^2}{3b}\,t,t),
\end{equation}
satisfies the complex modified Korteweg-de Vries type equation
\begin{equation}\label{1.21}
\begin{cases}
\p_t v+b\,\p_x^3 v+d\,|v|^2\p_x v+ e\,v^2\p_x \bar v=0,\\
v(x,0)  = \exp(i a\,x/3b)\,
u(x,0).
\end{cases}
\end{equation}
%If and only if $u$ solves (\ref{1.1})
Was shown in
\cite{C2} that the flow associated to the IVP (\ref{1.1})
leaves the following quantity
\begin{equation}
I_1(u) = \int_{\R}|u|^2(x,t)\,dx, \label{1.2}
\end{equation}
conserved in time. Also, when $be \neq 0$ we
have the following conserved quantity
\begin{align}
I_2(u) = c_1\int_{\R}|\partial_x u|^{2}(x,t)\,dx
+c_2\int_{\R}|u|^4(x,t)dx + c_3\textrm{Im}\int_{\R}u(x,t)\partial_x\overline{u(x,t)}dx, \,
\label{1.3}
\end{align}
where $c_1=3be$, $c_2=-e(e+d)/2$ and $c_3=(3bc-a(e+d))$.
We may suppose $c_3 = 0$. In fact, when $c_3 \neq 0$ we can take in the gauge transformation (\ref{GT})
$$\lambda= -\frac{c_3}{6be}.$$
Then, $u$ solves (\ref{1.1}) if and only if $v$
satisfies (\ref{1.9}) and in this new IVP we have the constant $c_3=0$.
  
We say that the IVP (\ref{1.1}) is locally well-posed in $X$
(Banach space) if the solution uniquely exists in certain time
interval $[-T,T]$ (unique existence), the solution describes a
continuous curve in $X$ in the interval $[-T,T]$ whenever initial
data belongs to $X$ (persistence), and the solution varies
continuously depending upon the initial data (continuous
dependence) i.e. continuity of application $u_{0} \mapsto u(t)$
from $X$ to $\mathcal{C}([-T,T];X)$.  We say that the IVP
(\ref{1.1}) is globally well-posed in $X$ if the same properties
hold for all time $T>0$. If some hypotheses in the definition of
local well-posed fail, we say that the IVP is ill-posed.

Particular cases of (\ref{1.1}) are the following:

$\bullet$
 Cubic nonlinear  Schr\"odinger equation
(NLS), ($a=\mp 1$, $b=0$, $c=-1$, $d=e=0$).
\begin{align}\label{3y0}
iu_{t} \pm u_{xx} + |u|^{2}u =  0, \quad x,t \in \R.
\end{align}
The best known local result for the IVP associated to (\ref{3y0})
is in $H^{s}(\R)$, $s \geq 0$, obtained by Tsutsumi \cite{Ts}.
Since the $L^2$ norm is preserved in (\ref{3y0}), one has that (\ref{3y0}) is globally well-posed in $H^{s}(\R)$, $s \geq 0$.

$\bullet$ Nonlinear Schr\"odinger equation with derivative
($a=-1 $, $b=0$, $c=0$, $d=2e$).
\begin{align}\label{4y0}
iu_{t} + u_{xx} + i\lambda(|u|^{2}u)_{x}=0, \quad x,t\in \R.
\end{align}
The best known local result for the IVP associated to  (\ref{4y0})
is in $H^{s}(\R)$, $s \geq 1/2$, obtained by Takaoka \cite{T1}.
Colliander et al. \cite{C-K-S-T-T3} they proved that (\ref{4y0}) is globally well-posed in $H^{s}(\R)$, $s >1/2$.

$\bullet$ Complex modified  Korteweg-de Vries (mKdV)
equation  ($a=0$, $b=1$, $c=0$, $d=1$, $e=0$).
\begin{align}\label{83y0}
u_{t}+u_{xxx}+|u|^{2}u_{x}=0, \quad x,t \in \R.
\end{align}
If $u$ is real, (\ref{83y0}) is the usual mKdV equation. Kenig
et al. \cite{KPV1} proved that the IVP associated to it is locally
well-posed in $H^{s}(\R)$, $s \geq 1/4$ and Colliander et al. \cite{C-K-S-T-T2}, proved that \eqref{83y0} is globally well-posed in $H^{s}(\R)$, $s >1/4$.

$\bullet$ When $a \neq 0$ is real and $b=0$, we
obtain a particular case of the well-known mixed nonlinear
Schr\"odinger equation
\begin{align}\label{6y0}
u_{t}=i a u_{xx}+\lambda(|u|^{2})_{x}u+g(u),\quad x,t \in \R,
\end{align}
where $g$ satisfies some appropriate conditions and $\lambda \in \R$ is a constant. Ozawa and
Tsutsumi in \cite{[O-T]} proved that for any $\rho >0$, there is a
positive constant $T(\rho)$ depending only on $\rho$ and $g$, such
that the IVP (\ref{6y0}) is locally well-posed in $H^{1/2}(\R)$,
whenever the initial data satisfies
\[\|u_{0}\|_{\mathrm{H}^{1/2}} \leq \rho.\]

There are other dispersive models similar to (\ref{1.1}),
see for instance \cite{[A],[CC],[P-Na],[P-S-S-M],[SS]}
and the references therein.

Regarding the IVP (\ref{1.1}), Laurey in \cite{C2} showed that
the IVP is locally well-posed in $H^{s}(\R)$ with $s > 3/4$, and using the quantities \eqref{1.2} and \eqref{1.3} she proved the global well-posedness
in $H^s(\R)$ with $s\geq 1$.
In \cite{G} Staffilani established the local well-posedness in $H^{s}(\R)$ with $s\geq 1/4$, for the IVP associated to (\ref{1.1}),
improving Laurey's result.

In the IVP (\ref{1.1}), when $a$, $b$ are real functions of $t$,
in \cite{Cv1, CL} was prove the local well-posedness in
$H^{s}(\R)$, $s \geq 1/4$.
Also, in \cite{Cv1, CM} was study the unique continuation property
for the solution of (\ref{1.1}).

\begin{remark}
{\bf 1)} Using \eqref{1.20} and the results obtained in \cite{C-K-S-T-T2}
we have that the PVI \eqref{1.1} is globally well-posed in
$H^s(\R)$ with $s> 1/4$, for initial data of the form:
$$\exp{\{\textstyle -i\frac{\textstyle a}{\textstyle 3b}x\}}v_0(x), \quad \exp{\{\textstyle -i\frac{\textstyle a}{\textstyle 3b}x\}}(v_0(x)+iv_0(x)),$$
where $v_0 \in H^s$, $s>1/4$, $v_0 \in \R$.
Therefore it suggests us to improve the result and obtain the global existence
for the general case in $H^s(\R)$, $s> 1/4$.

{\bf 2)}
If \, $e=0$, \,$b  d>0$ and 
$c=(a/3b)d$ in (\ref{1.1}), then the equation 
\begin{eqnarray}\label{11y0} 
  \partial_{t}u+ia
\partial^{2}_{x}u+b \partial^{3}_{x}u+i \frac{a}{3b}d|u|^{2}u+d
|u|^{2}\partial_{x}u  =  0,
\end{eqnarray}
have the following solution with two parameters 
\begin{align}\label{12y0}
u_{\eta,N}(x,t)= f_{\eta}(x+\psi(\eta,N)t)\exp i\{Nx+\phi(\eta,N)t\},
\end{align}
where $f_{\eta}(x)=\eta f(\eta x)$, $f(x)= (A \cosh{x})^{-1}$, $A=\sqrt{d/(6b)}$, $
\psi(\eta,N)=2a N+3b
N^{2}-\eta^{2}b$ and  $\phi(\eta,N)=aN^{2}+b N^{3}-3\eta^{2}
b N -a \eta^{2}$.

Using the transformation (\ref{1.20}) we can to obtain other family of solutions for (\ref{11y0}). In
fact, let $w$ solution of
\begin{eqnarray}\label{y0}
\left\{ \begin{array}{ll} \partial_{t}w+  
\partial^{3}_{x}w + 
|w|^{2}\partial_{x}w = 0, \quad x,t \in \R, &  \\ 
w(x,0)  =  w_{0}(x)=f_{1}(x)\exp i\{Nx\}=(\frac{1}{\sqrt{6}} \cosh{x})^{-1}\exp i\{Nx\},  
& \end{array} \right.
\end{eqnarray}
given by (\ref{12y0}). If $w$ is a solution of (\ref{y0}), then
$$v(x,t)=\frac{1}{\alpha} w(b^{-1/3}x,t), \quad \alpha=\sqrt{\frac{d}{b^{1/3}}}$$ is a solution of
\begin{eqnarray}\label{y1}
\left\{ \begin{array}{ll} \partial_{t}v+  
b\partial^{3}_{x}v + d
|v|^{2}\partial_{x}v= 0, \quad x,t \in \R, &  \\ 
v(x,0)  =  v_{0}(x),  
& \end{array} \right.
\end{eqnarray}
with initial data $v_{0}(x)=(1/\alpha)w(b^{-1/3}x,0)$ and if $v$ is a solution of (\ref{y1}) then,
using the transformation (\ref{1.20}) 
$$u(x,t)=v(x-\frac{a^2}{3b}t,t)\exp{i(\frac{2a^3}{27\b^2}t-\frac{a}{3b}x) } $$
is a solution of (\ref{11y0})
with initial data $u_{0}(x)=v(x,0)\exp{\{-i(a/3b)x\}}$,
therefore other solution of (\ref{11y0}) with two parameters is
\begin{align}
u_{\eta,N}(x,t)= g_{\eta}(b^{-1/3}x+\psi(\eta,N)t)\exp{\{ix(b^{-1/3}N-\frac{a}{3b})+it
\phi(\eta,N)\}},
\end{align}
where $g(x)=(\tilde{\alpha}\cosh{x})^{-1}$, $\tilde{\alpha}=\alpha/\sqrt{6}$, $\phi(\eta,N)=2a^3/(27b^2)-3N \eta^2+N^3 -Na^2 b^{-1/3}/(3b)$,
$\psi(\eta,N)=-a^2 b^{-1/3}/(3b)-\eta^2+3N^2$ and
$$u_{\eta,N}(x,0)  =  {u_{0}}_{\eta,N}(x)= g_{\eta}(b^{-1/3}x)\exp{\{ix(b^{-1/3}N-\frac{a}{3b})\}}.$$

When $a=0$ and $b=d=1$ in (\ref{11y0}), this solution coincide with the solution obtained in \cite{KPV5}.

{\bf 3)} If \;$e \neq 0$ and $b(d+e)>0$, then (\ref{1.1}) have solutions with one parameter:
\begin{align*}
u_{\eta}(x,t)= g_{\eta}(x+\psi(\eta,w)t)\exp i\{wx+\phi(\eta,w)t\},
\end{align*}
where
 $w= (c-2a A^{2})/(2e)$, $g_{\eta}(x)=\eta g(\eta x)$,
$g(x)=(A\cosh{x})^{-1}$, $A= \sqrt{(e +d)/(6b)}$, $\psi$ and $\phi$ as in (\ref{12y0}).

We have also that if $u$ is a solution of (\ref{1.1}) then, $v=\alpha u$ is a solution of (\ref{1.1}), where $\alpha \in
\C$, $|\alpha|=1$,  and if
$d \neq e$ in (\ref{1.1}) then $u(x,t)=\exp{i\{Cx+Dt+C_{0}) \}}$ is a solution of (\ref{1.1}),
where $D=\a C^{2}+b C^{3}$ e $C=c/(e-d)$.
\end{remark}
Recently there appeared several papers devoted to the global solution
of the dispersive type equation, where the framework is based
on almost conserved laws and the I-method,
see \cite{C-K-S-T-T1,C-K-S-T-T3,C-K-S-T-T2,FLP,FLP1}.
In this paper we adopt this way in order to obtain our results.

Our aim in this paper is to extend the local solution to a global one. Now, we state our main theorem of global existence:
\begin{theorem}\label{Tglobal}
The IVP (\ref{1.1}), with $c=(d-e)a/3b$, is global well-posedness in $H^s$, $s>1/4$.
\end{theorem}
{\bf{Notation.}} The notation to be used is mostly standard.
We will use the space-time Lebesgue $L_{x}^{p}L_{T}^{q}$ endowed
with the norm
$$
\|f\|_{L_{x}^{p}L_{T}^{q}} =
\big\| \|f\|_{L_{T}^{q}} \big\|_{L_{x}^{p}}
= \Big( \int_{\R} \Big( \int _{0}^{T} |f(x,t)|^{q} dt \Big)^{p/q} dx \Big)^{1/p}.
$$
We will use the notation $\|f\|_{L_{x}^{p}L_{t}^{q}}$ when the integration in the time variable is on the whole real line. In order to define the $X_{s, \beta}$ spaces we consider the following
IVP
\begin{equation*}%\label{eq:hs0}
  \left \{
   \begin{array}{l}
     u_t+i au_{xx}+bu_{xxx}=0, \quad x, \,t\in \mathbb{R},\,\,b \neq 0,\\
     u(0)=u_0,
   \end{array}
  \right.
\end{equation*}
whose solution is given by $u(x,t)=U(t)u_0(x)$, where the unitary
group $U(t)$ is defined as
\begin{align*}%\label{gU}
\widehat{U(t)u_0}(\xi)=e^{it(b\xi^3+a\xi^2)}\widehat{u_0}(\xi).
\end{align*}
For $s,\beta\in \mathbb{R}$, $X_{s,\beta}$ denotes the completion of the
Schwartz space $S(\mathbb{R}^2)$ with respect to the norm
\begin{align*}%\label{eq:norm1}
  \|u\|_{s,\beta}\equiv \|u\|_{X_{s,\beta}} \equiv  \|U(-t)u\|_{H_{s,\beta}} \equiv & \|\langle \tau \rangle^{\beta}
  \langle \xi \rangle^{s}  \widehat{U(-t)u}(\xi,\tau) \|_{L_{\tau}^2L_{\xi}^2} \\
  = & \|\langle \tau-(b\xi^3+a\xi^2) \rangle^{\beta}
  \langle \xi \rangle^{s}  \widehat u(\xi,\tau) \|_{L_{\tau}^2L_{\xi}^2 {\textstyle,}}
\end{align*}
where $$\widehat u(\xi,\tau)\equiv \int_{\mathbb{R}^2}
e^{-i(x\xi+t\tau)}u(x,t)dxdt.$$
For any time interval $[0,\rho]$, we define the space $X_{s,b}^{\rho}$ by the norm
\begin{equation*}
\|u\|_{X_{s,b}^{\rho}}=\inf\{\|U\|_{X_{s,b}}:\, U|_{[0,\rho]\times \R}=u\}.
\end{equation*}
The notation $A \lesssim B$ means there exist a constant $C$ such that $A\leq C \;B$,
and $A \thicksim B$ means $A\lesssim B$ and $B \lesssim A$. The notations $\xi_{ij}$ means $\xi_i+\xi_j$, $\xi_{ijk}$ means $\xi_i+\xi_j+\xi_k$, etc. Also we use the notation $m(\xi_i):=m_i$, $m(\xi_{ij}):=m_{ij}$, etc.

The notations for multilinear expressions is the same as in \cite{C-K-S-T-T1, C-K-S-T-T3}, we define a spatial n-multiplier to be any function $M_n(\xi_1, \ldots, \xi_n)$ on the hyperplane 
\begin{align*}
\Gamma_n:=\{(\xi_1, \ldots, \xi_n) \in \R^n; \xi_1+ \cdots+\xi_n=0\},
\end{align*}
which we endow with the dirac measure $\delta(\xi_1+\cdots+\xi_n)$. 
We define the n-linear functional  as
\begin{align*}
\Lambda(M_n;f_1, \ldots,f_n):=\int_{\Gamma_n}M_n(\xi_1, \ldots, \xi_n)\prod_{1}^{n}\widehat{f}_j(\xi_j),
\end{align*}
where $f_1, \ldots,f_n$ are complex functions on $\R$. We shall denote
\begin{align*}
\Lambda(M_n;f):=\Lambda(M_n;f,\overline{f},f,\overline{f} \ldots,f,\overline{f}).
\end{align*}
For $1\le j \le n$, $k \ge 1$ we define the elongation $\textbf{X}_j^k(M_n)$ of $M_n$ to be the multiplier of order $n+k$ given by
\begin{align*}
\textbf{X}_j^k(M_n)(\xi_1, \ldots, \xi_{n+k}):=M_n(\xi_1, \ldots, \xi_{j-1},\xi_j, \ldots, \xi_{j+k},\xi_{j+k+1}, \ldots, \xi_{n+k}).
\end{align*}

%%%%%%%%%%%%%%%%%%%%%%%%%%%%%%%%%%%%%%%%%%%%%%%%%%%%%%%%%%%%%%%%
\secao{Almost Conservations Laws}
%%%%%%%%%%%%%%%%%%%%%%%%%%%%%%%%%%%%%%%%%%%%%%%%%%%%%%%%%%%%%%%%
From (\ref{1.1}) we have
\begin{align*}
\p_t w+ia\,\p^2_x w+b\,\p^3_x w+ic\,w\bar{w}w+d\,
 (\p_x w) \bar{w}w+e\, w(\p_x\bar{w})w=&  0,\\
\p_t \bar{w}-ia\,\p^2_x \bar{w}+b\,\p^3_x \bar{w}-ic\,\bar{w}w\bar{w}+d\,
 (\p_x \bar{w})w \bar{w}+e\, \bar{w}(\p_x w)\bar{w}=&  0.
\end{align*}
Taking Fourier transformation in the above equalities we obtain the following result
\begin{proposition}
Let $n \ge 2$, be an even integer, and let $M_n$ be a multiplier of order $n$, then
\begin{align}\label{ame}
\p_t \Lambda_n(M_n;w)=&i \Lambda_n(M_n \Upsilon_n^{a,b};w)-i\Lambda_{n+2}
\left(\sum_{j=1}^{n}\Upsilon_{j,n+2}^{c,e}\textbf{X}_j^2(M_n);w\right)\nonumber\\
&-i d \Lambda_{n+2}\left(\sum_{j=1}^{n/2}\textbf{X}_{2j-1}^2(M_n)\xi_{2j-1}+\sum_{j=1}^{n/2}\textbf{X}_{2j}^2(M_n)\xi_{2j+2};w\right),
\end{align} 
where $\Upsilon_n^{a,b}=\sum_{j=1}^{n}((-1)^{j-1}a\xi_{j}^2+ b\xi_{j}^3)$ and $\Upsilon_{j,n+2}^{c,e}=(-1)^{j-1}c+ e\xi_{j+1}$.
\end{proposition} 
We define the first modified energy as
\begin{align}\label{E1}
E_1=k_1 \Lambda_2(M_2;w), \quad M_2(\xi_1,\xi_2)=\xi_1 \xi_2 m(\xi_1)m(\xi_2),
\end{align}
where $k_1=3be$, and the second modified energy as
\begin{align}\label{E2}
E_2=E_1+ \Lambda_4(\delta_4),
\end{align}
where the 4-multiplier $\delta_4$ will be choosed after. By (\ref{ame}) we get
\begin{align}
\p_t E_2=& \p_t E_1+ \p_t \Lambda_4(\delta_4)= k_1 \Lambda_2(M_2\Upsilon_2^{a,b})-ik_1\Lambda_{4}
\left(\sum_{j=1}^{2}\Upsilon_{j,4}^{c,e}\textbf{X}_j^2(M_2)\right) \nonumber\\
& -id k_1\Lambda_{4}(\textbf{X}_{1}^2(M_2)\xi_{1}+\textbf{X}_{2}^2(M_2)\xi_{4})+i\Lambda_{4}(\delta_4\Upsilon_4^{a,b})-i\Lambda_{6}
\left(\sum_{j=1}^{4}\Upsilon_{j,6}^{c,e}\textbf{X}_j^2(\delta_4)\right)\nonumber\\
& -id \Lambda_{6}\left(\sum_{j=1}^{2}\textbf{X}_{2j-1}^2(\delta_4)\xi_{2j-1}+\sum_{j=1}^{2}\textbf{X}_{2j}^2(\delta_4)\xi_{2j+2}\right),\label{Et}
\end{align}
it is clear that $\Lambda_2(M_2\Upsilon_2^{a,b})=0$. Now if $\tilde{M_n}$ is an n- multiplier 
$\Lambda_n(\tilde{M_n})$ is invariant under permutations of the even $\xi_j$ indices or of the odd $\xi_j$ 
indices, therefore for achieve a cancellation of the 4-linear expression, we choose $\delta_4$ such that
\begin{align}
\Upsilon_4^{a,b} \delta_4=& \frac{ck_1}{2} (\xi_1^2m_1^2-\xi_2^2m_2^2+\xi_3^2m_3^2-\xi_4^2m_4^2)-\frac{ek_1}{2}(\xi_2 \xi_4^2 m_4^2+\xi_4 \xi_2^2 m_2^2+\xi_1 \xi_3^2 m_3^2+\xi_3 \xi_1^2 m_1^2)\nonumber\\
& -\frac{dk_1}{2}(\xi_1 \xi_4^2 m_4^2+\xi_4 \xi_1^2 m_1^2+\xi_3 \xi_2^2 m_2^2+\xi_2 \xi_3^2 m_3^2),\label{X4}
\end{align} 
consequently from (\ref{Et}) we get
\begin{align}\label{d6}
\p_t E_2=\Lambda_6(\delta_6),
\end{align}
with

\begin{align}
\delta_6=& \frac{-ie}{36}	
\mathop{\sum\limits_{\{k,m,o\}=\{1,3,5\}}}_{\{l,n,p\}=\{2,4,6\}}
%\{k,m,o\}=\{1,3,5\} ; \{l,n,p\}=\{2,4,6\}
[
\xi_{l}\delta_4(\xi_{klm},\xi_{n},\xi_{o},\xi_{p} )+\xi_m \delta_4(\xi_{k},\xi_{lmn},\xi_{o},\xi_{p} )+\xi_n \delta_4(\xi_{k},\xi_{l},\xi_{mno},\xi_{p} )\nonumber\\
&+\xi_o \delta_4(\xi_{k},\xi_{l},\xi_{m},\xi_{nop} )
]
-\frac{id}{36}	
\mathop{\sum\limits_{\{k,m,o\}=\{1,3,5\}}}_{\{l,n,p\}=\{2,4,6\}}
%\{k,m,o\}=\{1,3,5\} ; \{l,n,p\}=\{2,4,6\}
[
\xi_{k}\delta_4(\xi_{klm},\xi_{n},\xi_{o},\xi_{p} )+\xi_m \delta_4(\xi_{k},\xi_{l},\xi_{mno},\xi_{p})\nonumber\\ 
&+\xi_n \delta_4(\xi_{k},\xi_{lmn},\xi_{o},\xi_{p})+\xi_p \delta_4(\xi_{k},\xi_{l},\xi_{m},\xi_{nop})
].
\nonumber
\end{align}
\begin{proposition}\label{m1}
If $m(\xi)=1$ for all $\xi$, then
$$\p_t E_2=0.$$
\end{proposition}
\begin{proof}
From definition of $E_2$, we have
\begin{align}\label{E21}
E_2= 3be \Lambda_2(\xi_1 \xi_2;w)+ \Lambda_4(\delta_4;w),
\end{align}
where
\begin{align*}
\Upsilon_4^{a,b} \delta_4=& \frac{ck_1}{2} (\xi_1^2-\xi_2^2+\xi_3^2-\xi_4^2)-\frac{ek_1}{2}(\xi_2 \xi_4^2 +\xi_4 \xi_2^2 +\xi_1 \xi_3^2 +\xi_3 \xi_1^2 )\\
& -\frac{dk_1}{2}(\xi_1 \xi_4^2 +\xi_4 \xi_1^2+\xi_3 \xi_2^2 +\xi_2 \xi_3^2 ).
\end{align*} 
If $\xi_1+\cdots +\xi_4=0$, then $\xi_1^2-\xi_2^2+\xi_3^2-\xi_4^2=2\xi_{12}\xi_{14}$ and $\xi_1^3+\cdots +\xi_4^3= 3\xi_{12}\xi_{13}\xi_{14}$, therefore
\begin{align}\label{upsilon}
\Upsilon_4^{a,b}=2a\xi_{12}\xi_{14}+3b\xi_{12}\xi_{13}\xi_{14}=\xi_{12}\xi_{14} (2a+3b\xi_{13}).
\end{align}
On the other hand $\xi_2 \xi_4^2 +\xi_4 \xi_2^2 +\xi_1 \xi_3^2 +\xi_3 \xi_1^2 =-\xi_{12}\xi_{13}\xi_{14}$\, (see Lemma 3.5 in \cite{C-K-S-T-T1} and Remark 3.6 in \cite{C-K-S-T-T3}), similarly $\xi_1 \xi_4^2 +\xi_4 \xi_1^2+\xi_3 \xi_2^2 +\xi_2 \xi_3^2=-\xi_{12}\xi_{13}\xi_{14}$, hence
\begin{align*}
\delta_4= &\frac{3be}{2}\frac{2c \xi_{12}\xi_{14} +(d+e)\xi_{12}\xi_{13}\xi_{14}}{\xi_{12}\xi_{14} (2a+3b\xi_{13})}\\
=&\frac{e(d+e)}{2}.
\end{align*}
And from (\ref{E21}) we get
\begin{align*}
E_2=& 3be \Lambda_2(\xi_1 \xi_2;w)+ \frac{e(d+e)}{2}\Lambda_4(1;w)\\
=&-3be\int_{\R}|w_x|^2+ \frac{e(d+e)}{2}\int_{\R}|w|^4\\
=&-I_2(w).
\end{align*}
This concludes the proof of the proposition.
\end{proof}
In the following sections we will consider $a=c=0$ in the IVP (\ref{1.1}) (see (\ref{1.20}) and (\ref{1.21})).
%%%%%%%%%%%%%%%%%%%%%%%%%%%%%%%%%%%%%%%%%%%%%%%%%%%%%%%%%%%%%%%%
\secao{Preliminary results}
%%%%%%%%%%%%%%%%%%%%%%%%%%%%%%%%%%%%%%%%%%%%%%%%%%%%%%%%%%%%%%%%
For the estimates on the multipliers we use the following elementary results.  
\begin{lemma}\label{TMVT}{\bf 1)} (Double mean value theorem DMVT)\\
Let $f \in C^2(\R)$, and $\max{\{|\eta|, |\lambda|\}}\ll \xi $, then
\begin{align*}
|f(\xi+\eta+\lambda)-f(\xi+\eta)-f(\xi+\lambda)+f(\xi)|\lesssim |f''(\theta)||\eta||\lambda|, 
\end{align*}
where $|\theta| \sim |\xi|$.

{\bf 2)} (Triple mean value theorem TMVT) \\Let $f \in C^3(\R)$, and $\max{\{|\eta|, |\lambda|, |\gamma|\}}\ll \xi $, then
\begin{equation*}
|f(\xi+\eta+\lambda+\gamma) - f(\xi+\lambda+\eta)-f(\xi+\eta +\gamma)-f(\xi+\lambda+\gamma)+f(\xi+\eta)
\end{equation*}
\begin{equation*}+f(\xi +\lambda)+f(\xi+\gamma)-f(\xi)|\lesssim |f'''(\theta)||\eta||\lambda||\gamma|,
\end{equation*}
where $|\theta| \sim |\xi|$.
\end{lemma}
And for the proof of Proposition \ref{lemiter}, shall be fundamental the improved Strichartz estimate.
\begin{lemma}\label{4-2}
Let $s>1/4$, $v_1,v_2 \in \mathbf{S}(\R\times\R)$ such that $\textrm{supp}\,{\widehat{v_1}} \subset \{|\xi| \sim N\}$ and $\textrm{supp}\,{\widehat{v_2}} \subset \{|\xi| \ll N\}$, then
\begin{align*}
\|v_1 v_2\|_{L_x^4L_t^2} \lesssim \frac{1}{(1-4s)^{1/4}}\frac{1}{N} \|v_2\|_{X_{s,1/2+}^{\rho}}\|v_1\|_{X_{0,1/2+}^{\rho}}.
\end{align*}
\end{lemma}
\begin{proof}
As in \cite{C-K-S-T-T1} is sufficient to prove 
\begin{align*}
\|v_1 v_2\|_{L_x^4L_t^2} \lesssim \frac{1}{(1-4s)^{1/4}}\frac{1}{N} \|\phi\|_{H^{s}}\|\psi\|_{L^2},
\end{align*}
where $v_1=U(t)\psi$ and $v_2=U(t)\phi$ . By duality, definition of $v_1, v_2$, Fubini theorem and Plancherel identity in the spatial variable we have
\begin{align*}
\|v_1 v_2\|_{L_x^4L_t^2} \lesssim &  \sup_{\|F\|_{L_x^{4/3}L_t^2}\le 1}\int_{\R^2}|\widehat{\phi}(y)\widehat{\psi}(z)\widehat{F}(z+y,z^3+y^3)| dz dy\\
=& \frac{1}{N^{2}} \sup_{\|F\|_{L_x^{4/3}L_t^2}\le 1}\int_{\R^2}|\widehat{\phi}(y(s,r))\widehat{\psi}(z(s,r))\widehat{F}(s,r)| dr ds,
\end{align*}
where we used the change of variable $z+y=s$, $z^3+y^3=r$, which has Jacobian of size $N^2$.
Now if we applying H\"older inequality and a change of variables back for $z$ and $y$, we obtain
\begin{align*}
\|v_1 v_2\|_{L_x^4L_t^2} \lesssim & \frac{1}{N^2} \sup_{\|F\|_{L_x^{4/3}L_t^2}\le 1}\|\widehat{\phi}(y(s,r))\widehat{\psi}(z(s,r))\|_{L_s^{4/3}L_r^2} \|\widehat{F}\|_{L_x^{4}L_t^2}\\
\le & \frac{1}{N} \sup_{\|F\|_{L_x^{4/3}L_t^2}\le 1}\|\psi\|_{L_z^2} \|\widehat{\phi}\|_{L_y^{4/3}} \|\widehat{F}\|_{L_x^{4}L_t^2},
\end{align*}
where the Fourier transform of $F$ is taking only in the space variable. Using H\"older inequality we obtain for $s>1/4$
\begin{align*}
\int_{\R}|\widehat{\phi}|^{4/3}\le 
\Big( \int_{\R}\langle \xi^2 \rangle^s |\widehat{\phi}|^{2}\Big)^{2/3}
\Big( \int_{\R}\frac{1}{\langle \xi^2 \rangle^{2s}} \Big)^{1/3},
\end{align*}
therefore
\begin{align*}
\|\widehat{\phi}\|_{L_y^{4/3}}\le \frac{1}{(1-4s)^{1/4}}\|\phi\|_{H^s},
\end{align*}
and by Hausdorff-Young inequality and Minkowsky integral inequality, we get
\begin{align*}
\|\widehat{F}\|_{L_x^{4}L_t^2}\le \|\widehat{F}\|_{L_t^2 L_x^{4}}\le \|F\|_{L_t^2 L_x^{4/3}}\le \|F\|_{L_x^{4/3}L_t^2}\le 1.
\end{align*}
This completes the proof.
\end{proof}

We define the Fourier multiplier operator I with symbol
\begin{align}\label{m33}
m(\xi)= 
\left\{\begin{array}{ll}
1, &  |\xi|<N,\\
\frac{\textstyle N^{1-s}}{\textstyle |\xi|^{1-s}}, & |\xi|>2N.
\end{array} \right.
\end{align}
We have $I: H^s \mapsto H^1$.
For the local result we define the Fourier multiplier operator L, with symbol
\begin{displaymath}
l(\xi)=m(\xi)\langle \xi \rangle^{1-s}= 
\left\{\begin{array}{ll}
\langle \xi \rangle^{1-s}, &  |\xi|<N,\\
\langle \xi \rangle^{1-s}\frac{\textstyle N^{1-s}}{\textstyle |\xi|^{1-s}}, & |\xi|>2N.
\end{array} \right.
\end{displaymath}
Is obvious that \begin{align}\label{opL}
\|Iu\|_{H^1}=\|Lu\|_{H^s},\quad  \|Iu\|_{X_{1,b}}=\|Lu\|_{X_{s,b}\textrm{,}}\end{align} and
for $s \in [0,1)$ is
 $1 \leq l(\xi) \lesssim N^{1-s}$, therefore
\begin{align*}%\label{Il}
\|u\|_{s',b'}\lesssim \|Iu\|_{s'-s+1,b'} \lesssim N^{1-s}\|u\|_{s',b'},\quad s\in [0,1),
\end{align*}
observe that if $V|_{[0,\rho]\times \R}=Iu$, $V \in X_{s'-s+1,b'}$, then $U$ defined by $\widehat{U}=(1/m)\widehat{V}$, satisfies
$$ \|U\|_{s',b'}\lesssim\|V\|_{s'-s+1,b'},$$
moreover in $[0, \rho]$ is $U|_{[0,\rho]\times \R}=u$, therefore
%\begin{align}
%\|u\|_{X_{s',b'}^{\rho}} \leq \|U\|_{X_{s',b'}}\lesssim %\|V\|_{s'-s+1,b'}
%\end{align}
\begin{align}\label{Il}
\|u\|_{X_{s',b'}^{\rho}}\lesssim \|Iu\|_{X_{s'-s+1,b'}^{\rho}}. %\lesssim N^{1-s}\|u\|_{s',b'},\quad s\in [0,1),
\end{align}
Also we have 
 \begin{align}\label{opL1}
l(\xi_1+\xi_2) \lesssim l(\xi_1)+l(\xi_2).
\end{align}
In fact, for see this, without lost of generality we can assume $|\xi_1|\ge |\xi_2|$, we consider two cases: \\
i) If $|\xi_1|\le N$, then we have $|\xi_1+\xi_2| \le 2N$, this implies $$l(\xi_1+\xi_2) \sim \langle \xi_1+\xi_2 
\rangle^{1-s}\le \langle \xi_1 \rangle^{1-s}+\langle \xi_2 \rangle^{1-s}=l(\xi_1)+ l(\xi_2).$$
ii) If $|\xi_1|\ge N$, then we have $l(\xi_1) \sim N^{1-s}$, thus for all $\xi$, $l(\xi) \lesssim l(\xi_1)$, 
in particular $l(\xi_1+\xi_2) \lesssim l(\xi_1) \le l(\xi_1)+l(\xi_2)$. Note that (\ref{opL1}) implies $l(\xi_1+\xi_2) \lesssim l(\xi_1)l(\xi_2).$

In the proof of Theorem \ref{Tglobal} we will use the following local result.
\begin{theorem}\label{tlocal}
Let $s \geq 1/4$, then the IVP (\ref{1.1}) is locally well-posed for data $\varphi$, with $I \varphi \in H^1$ where the time of existence satisfies
\begin{align}\label{3tempo3}
\delta \sim \|I \varphi\|_{H^1}^{-\theta},
\end{align}
with $\theta>0$. Moreover the solution of the IVP (\ref{1.1}), is such that
\begin{align}\label{3Iu3}
\|I u\|_{X_{1,1/2+}^{\delta}} \lesssim  \|I u_0\|_{H^1 \textrm{.}}
\end{align}
\end{theorem}
\begin{proof}
The Theorem \ref{tlocal} is practically done in \cite{Tk} (see also \cite{ZB}), in fact, is sufficient to prove
\begin{align}
\|L (uv\overline{w}_x)\|_{X_{s,-1/2+}}\lesssim & \;\|L u\|_{X_{s,1/2+}}\|L v\|_{X_{s,1/2+}}\|L w\|_{X_{s,1/2+}},\label{tloc1}\\
\|L (u\overline{v}w_x)\|_{X_{s,-1/2+}}\lesssim & \; |L u\|_{X_{s,1/2+}}\|I v\|_{X_{s,1/2+}}\|L w\|_{X_{s,1/2+}},\nonumber\\
\|L (u\overline{v}w)\|_{X_{s,-1/2+}}\lesssim & \; \|L u\|_{X_{s,1/2+}}\|L v\|_{X_{s,1/2+}}\|L w\|_{X_{s,1/2+}},\label{tloc2}
\end{align}
in order to prove the first inequality we make the following decomposition
\begin{align*}
	l(\xi)\widehat{uv\overline{w}_x}(\xi,\tau) = &l(\xi)\int\limits_{|\xi_1|>2N}\zeta+l(\xi)\mathop{\int\limits_{|\xi_1| \leq 2N}}_{|\xi_2|> 2N}\zeta+l(\xi)\mathop{\mathop{\int\limits_{|\xi_1|\leq 2N}}_{|\xi_2| \leq 2N}}_{|\xi-\xi_1-\xi_2| > 2N}\zeta\\
	&+l(\xi)\mathop{\mathop{\int\limits_{|\xi_1|\leq 2N}}_{|\xi_2| \leq 2N}}_{|\xi-\xi_1-\xi_2| \le 2N}\zeta,
\end{align*}
where $\zeta:=\widehat{u}(\xi-\xi_1-\xi_2,\tau-\tau_1-\tau_2)\widehat{v}(\xi_1,\tau_1)\widehat{\overline{w}}_x(\xi_2,\tau_2)$, thus
\begin{align*}
|l(\xi)\widehat{uv\overline{w}_x}(\xi,\tau)| \lesssim |\widehat{uv_1\widehat{\overline{w}}_x}|+ |\widehat{uv_2\widehat{\overline{v_3}}_x}| +|\widehat{v_2\widehat{\overline{v_4}}_xv_5}|+\langle \xi\rangle^{1-s} |\widehat{v_2\widehat{\overline{v_4}}_x v_6}|,
\end{align*}
with
\begin{align*}
&\widehat{v_1}(\xi,\tau)=\Chi_{|\xi|>2N}\widehat{v}(\xi,\tau)l(\xi),\quad \widehat{v_2}(\xi,\tau)=\Chi_{|\xi|\leq 2N}\widehat{v}(\xi,\tau),\\ &\widehat{\overline{v_3}}_x(\xi,\tau)=\Chi_{|\xi|>2N}\widehat{\overline{w}}_x(\xi,\tau)l(\xi), \quad
\widehat{\overline{v_4}}_x(\xi,\tau)=\Chi_{|\xi|\leq 2N}\widehat{\overline{w}}_x(\xi,\tau), \\
&\widehat{v_5}(\xi,\tau)=\Chi_{|\xi|> 2N}\widehat{u}(\xi,\tau)l(\xi),\quad \widehat{v_6}(\xi,\tau)=\Chi_{|\xi|\leq 2N}\widehat{u}(\xi,\tau), 
\end{align*}
 and applying Proposition 2.7 in \cite{Tk} (see also Theorem 2.1 in \cite{ZB}) we obtain (\ref{tloc1}).  For the other inequalities we make an analogous decomposition.
\end{proof}

We also have a result of local well-posed with the interval of existence as in (\ref{3tempo3}), without the use of the theory of the spaces $X_{s,b}$.

In fact, let 
\begin{align*}
|||u|||_{\Delta T, s}=&\|\p_{x}u\|_{L_x^{\infty}L_{\Delta\!T}^2}+
\|D_{x}^{s}\p_{x}u\|_{L_x^{\infty}L_{\Delta\!T}^2}+
\|D_{x}^{s-1/4}\p_x\,u\|_{L_x^{20}L_{\Delta\!T}^{5/2}} 
 +\|u\|_{L_x^5 L_{\Delta\!T}^{10}}\\
&+\|D_x^s u\|_{L_x^5 L_{\Delta\!T}^{10}} +
\|u\|_{L_x^4 L_{\Delta\!T}^{\infty}} 
 +\|u\|_{L_x^8L_{\Delta\!T}^8}
+ \|D_{x}^s \,u\|_{L_x^8 L_{\Delta\!T}^8\textrm{.}}
\end{align*}
\begin{theorem}\label{t3.1}
Let  $s \geq 1/4$,
and $a, b  \in \R, b\neq 0$, $c , d , e \in \C$, then the IVP (\ref{1.1}) is locally well-posed for data $\varphi$, with $I \varphi \in H^1$. Moreover the solution is such that
\begin{equation*}
|||u|||_{\delta,s} \lesssim \|I \varphi\|_{H^1 \textrm{,}}
\end{equation*}
where $\delta$ satisfies (\ref{3tempo3}) with $\theta=4$.
\end{theorem}
\begin{proof}
The theorem follows from the proof in \cite{G} if we prove
\begin{equation*}
|||U(t)u_0|||_{\delta,s} \lesssim \|I u_0\|_{H^1}=\|L u_0\|_{H^s \textrm{,}}
\end{equation*}
we consider the first term in the definition of $||| \cdot |||_{T}$, we will prove
\begin{align}\label{regloc}
\|\p_{x}U(t)u_0\|_{L_x^{\infty}L_{\delta}^2}\lesssim \|L u_0\|_{L^2}.
\end{align}
The inequality (\ref{regloc}) is equivalent with
\begin{align*}
\|L^{-1}\p_{x}U(t)u_0\|_{L_x^{\infty}L_{\delta}^2}\lesssim \| u_0\|_{L^2 \textrm{,}}
\end{align*}
where the Fourier multiplier operator $L^{-1}$ have symbol $1/l(\xi) \le 1$, it is easy to see that
\begin{align}
\|L^{-1}\p_{x}U(t)u_0\|_{L_x^{\infty}L_{\delta}^2}\le \|\p_{x}U(t)L^{-1}u_0\|_{L_x^{\infty}L_{\delta}^2} \le \|L^{-1} u_0\|_{L^2} \le \| u_0\|_{L^2 \textrm{.}}
\end{align}
We proceed similarly with the others terms.
\end{proof}
\begin{lemma}\label{strich}For any $s_1\geq 1/4$, $s_2 \geq 0$ and $b>1/2$ we have
\begin{align}
\|u\|_{L_x^{4}L_{\delta}^{\infty}} \lesssim & \, \|u\|_{X_{s_1, b}^{\rho}}, \label{4-1}\\
\|u\|_{L_x^{8}L_{\delta}^{8}} \lesssim & \, \|u\|_{X_{s_2, b}^{\rho}}, \label{8-8}\\
\|u\|_{L_x^{6}L_{\delta}^{6}} \lesssim & \, \|u\|_{X_{s_2, b}^{\rho}}. \label{6-6}
\end{align}
\end{lemma}
\begin{proof}
The inequalities (\ref{4-1}) and (\ref{8-8}), follows from 
\begin{align*}
\|U(t)u_0\|_{L_x^{4}L_{\delta}^{\infty}} \lesssim  \, \|u_0\|_{H^{1/4}}, \quad \|U(t)u_0\|_{L_x^{8}L_{\delta}^{8}} \lesssim  \, \|u_0\|_{L^{2}}
\end{align*}
and from a standard argument, see for example \cite{B2, C-K-S-T-T3, Gin}.

The inequality (\ref{6-6}) follows by interpolation between $\|v\|_{L_x^{8}L_{\delta}^{8}} \lesssim  \|v\|_{X_{0, 1/2+}}$ and the trivial estimate
$
\|v\|_{L_x^{2}L_{\delta}^{2}} \leq   \, \|v\|_{0,0}.
$
\end{proof}
%There differences between our method and the other's:\\
%i) The assumption on the low frequency part $v_0$ of $u_0$ (see \eqref{30})
%can be improved by using instead of
%$$
%\|v_0\|_{H^{\delta}}\leq \|u_0\|_{H^s}\,N^{\delta(1-s)},$$
%the estimate (\ref{31}). \\
%ii) Instead of Leibnitz's rule for fractional derivatives 
%(see \cite{FLP, FLP1}), we use the interpolation formula (\ref{2.5}).  \\
%iii) In \cite{Cv1}, in each iteration of size $\Delta T$,
%the problems (\ref{v}) and (\ref{w}) are defined in
%$[(j-1) \Delta T, j \Delta T]$ (j-th iteration).
%Now, those problems are defined in $[0,j \Delta T]$ for each iteration (see %the proof of Theorem \ref{t1.1}).
%In this way,  instead of the conservation law in $H^{1}$ (in our case (1.3))
%as in \cite{Cv2,CL,FLP}, we only make use of the Proposition \ref{limv}.
%\end{remark}
\begin{remark}
Actually the inequality (\ref{tloc2}) is valid for all $s>-1/4$ (see \cite{Cv2}).
\end{remark}
%%%%%%%%%%%%%%%%%%%%%%%%%%%%%%%%%%%%%%%%%%%%%%%%%%%%%%%%%%%%%%%%
\secao{Estimates for $\delta_4$ and $\delta_6$}
%%%%%%%%%%%%%%%%%%%%%%%%%%%%%%%%%%%%%%%%%%%%%%%%%%%%%%%%%%%%%%%%
From here onwards we will consider the notation
$|\xi_i|=N_i$, $m(N_i)=m_i$, $|\xi_{ij}|=N_{ij}$, $m(N_{ij})=m_{ij}$, etc. Given ford number $N_1$, $N_2$, $N_3$, $N_4$ and $\mathcal{C}=\{N_1, \ldots, N_4\},$ we will note
$N_s=\max \mathcal{C}$, $N_a=\max \mathcal{C} \setminus \{N_s\}$, $N_t=\max\mathcal{C}\setminus \{N_s, N_a\}$, $N_b=\min \mathcal{C}$, in this way $$N_s\ge N_a \ge N_t\ge N_b.$$

\begin{proposition}
Let $m$ defined as in (\ref{m33}), then
\begin{equation}\label{delta4}
|\delta_4| \lesssim m^2(N_s) 
\end{equation}
and
 \begin{equation}\label{delta6}
|\delta_6| \lesssim N_s m^2(N_s).
\end{equation}
\end{proposition}
In order to prove (\ref{delta4}) we will use the following proposition
in similar form like when $m=1$. 
\begin{proposition}\label{prop4.2} Let $m$ defined as in (\ref{m33}), then
\begin{equation}\label{4.3}
|\xi_2 \xi_4^2 m_4^2+\xi_4 \xi_2^2 m_2^2+\xi_1 \xi_3^2 m_3^2+\xi_3 \xi_1^2 m_1^2| \lesssim m^2(N_s)|\xi_{12}\xi_{13}\xi_{14}|,
\end{equation}
and
\begin{equation}\label{14.3}
|\xi_1 \xi_4^2 m_4^2+\xi_4 \xi_1^2 m_1^2+\xi_2 \xi_3^2 m_3^2+\xi_3 \xi_2^2 m_2^2| \lesssim m^2(N_s)|\xi_{12}\xi_{13}\xi_{14}|.
\end{equation}
\end{proposition}
\begin{proof}
Without lost of generality we can assume $|\xi_1|=N_s$, and by symmetry $|\xi_{12}| \le |\xi_{14}|$. In \cite{C-K-S-T-T3} (Lemma 4.1) they proved that
\begin{align*}
|\xi_2 \xi_4^2 m_4^2+\xi_4 \xi_2^2 m_2^2+\xi_1 \xi_3^2 m_3^2+\xi_3 \xi_1^2 m_1^2| \lesssim m^2(N_s)|\xi_{12}\xi_{14}|N_s, 
\end{align*} 
therefore we can suppose $|\xi_{13}|\ll N_s$, this implies $|\xi_{1}|\sim |\xi_{3}|$. Let $f(\xi)=\xi m(\xi)$,
observing that $\xi_{12}\xi_{14}=\xi_{2}\xi_{4}-\xi_{1}\xi_{3}$, we have
\begin{align}
\xi_2 \xi_4^2 m_4^2+\xi_4 \xi_2^2& m_2^2+\xi_1 \xi_3^2 m_3^2+\xi_3 \xi_1^2 m_1^2=\xi_2 \xi_4(f(\xi_2)+f(\xi_4))+\xi_1 \xi_3(f(\xi_1)+f(\xi_3))\nonumber\\
=&\xi_2\xi_4(f(\xi_1)+f(\xi_2)+f(\xi_3)+f(\xi_4))-\xi_{12}\xi_{14}(f(\xi_1)+f(\xi_3)).\label{m5}
\end{align}
In the second term of (\ref{m5}) we can use the medium value theorem (MVT) for to obtain
$$|\xi_{12}\xi_{14}(f(\xi_1)+f(\xi_3))|=|\xi_{12}\xi_{14}(f(\xi_1)-f(-\xi_3))| \lesssim |\xi_{12}\xi_{14}\xi_{13}|m^2 (N_s),$$ 
where we used that $|\xi_{13}| \ll N_s$, and $|f'(\xi_1)| \sim |m^2(\xi_1)|$. Therefore we will only estimate the first term in (\ref{m5}).

We consider two cases:
\\
1)\, $|\xi_{14}|\gtrsim |\xi_{3}|$, in this case we consider two sub-cases \\
a) If $|\xi_{12}|\ll |\xi_1| $, then using the DMVT (Lemma \ref{TMVT}) with $\xi=-\xi_1$, $\lambda=\xi_{12}$ and $\eta = \xi_{13}$ 
\begin{align*}
|\xi_2\xi_4(f(\xi_1)+f(\xi_2)+f(\xi_3)+f(\xi_4))| \lesssim |\xi_{14}|N_s|\xi_{12}\xi_{13}|\frac{m^2(N_s)}{N_s},
\end{align*}
where we also used that $|\xi_{2}| \le |\xi_{1} |\sim |\xi_{3} | \lesssim |\xi_{14}|$ and $|f''(\xi_1)| \lesssim m^2(\xi_1)/|\xi_1|$. \\
b) If $|\xi_{12}|\gtrsim |\xi_1|$, here we proceed similarly as in \cite{C-K-S-T-T3} (Lemma 4.1). Using the fact that $N_s \lesssim |\xi_{12}| \le |\xi_{14}|$, $(m^2(\xi) \xi^2)' \sim m^2(\xi) \xi$, $m^2(\xi) \xi$ is nondecreasing and the MVT we have
\begin{align*}
|\xi_2 \xi_4^2 m_4^2+\xi_4 \xi_2^2& m_2^2+\xi_1 \xi_3^2 m_3^2+\xi_3 \xi_1^2 m_1^2|=|\xi_3(m_1^2 \xi_1^2 -m_{1-13}^2\xi_{1-13}^2)+ \xi_2(m_4^2 \xi_4^2 \\
&-m_{4+13}^2\xi_{4+13}^2)-\xi_{24}(m_3^2 \xi_3^2 -m_{4+13}^2\xi_{4+13}^2)| \lesssim |\xi_{13}|N_s^2m^2(N_s).
\end{align*}
Hence we obtain (\ref{4.3}) in this sub-case.\\
2)\, $|\xi_{14}|\ll |\xi_{3}|$, using the TMVT considering in Lemma \ref{TMVT}: $\xi=-\xi_1$, $\lambda=\xi_{12}$, $\eta = \xi_{13}$ and $\gamma=\xi_{14}$ we have
\begin{align*}
|\xi_2\xi_4(f(\xi_1)+f(\xi_2)+f(\xi_3)+f(\xi_4))| \lesssim N_s^2 |\xi_{12}\xi_{13}\xi_{14}|\frac{m^2(N_s)}{N_s^2},
\end{align*}
where we also used that $|f'''(\xi_1)| \lesssim m^2(\xi_1)/|\xi_1|^2$.

Now, in order to obtain (\ref{14.3}), using (\ref{4.3}) we get
\begin{equation*}%\label{14.3}
|\xi_1 \xi_4^2 m_4^2+\xi_4 \xi_1^2 m_1^2+\xi_2 \xi_3^2 m_3^2+\xi_3 \xi_2^2 m_2^2| \lesssim m^2(N_s)|\xi_{21}\xi_{23}\xi_{24}|=m^2(N_s)|\xi_{12}\xi_{14}
\xi_{13}|.
\end{equation*}
This completes the proof.
\end{proof}
By (\ref{X4}), (\ref{upsilon}) and Proposition \ref{prop4.2}, we have (\ref{delta4}). The estimate (\ref{delta6}) is obvious.

%%%%%%%%%%%%%%%%%%%%%%%%%%%%%%%%%%%%%%%%%%%%%%%%%%%%%%%%%%%%%%%%
\secao{Estimates 4-lineal and 6-lineal}
%%%%%%%%%%%%%%%%%%%%%%%%%%%%%%%%%%%%%%%%%%%%%%%%%%%%%%%%%%%%%%%%
The following lemma will be used frequently in the estimates 4-lineal and 6-lineal.
\begin{lemma}
Let $n\ge 2$ a even integer, $w_1, \ldots, w_n \in \mathbf{S}(\R)$, then
\begin{equation}\label{int1}
\int_{\xi_1+\cdots +\xi_n=0}\widehat{w}_1\widehat{\overline{w}}_2\ldots\widehat{w}_{n-1}\widehat{\overline{w}}_n=\int_{\R}w_1\overline{w}_2\ldots w_{n-1}\overline{w}_n.
\end{equation}
\end{lemma}
In the proof of our global result, we will need the following properties.
\begin{proposition}\label{lemiter}
Let $w \in \mathbf{S}(\R \times \R)$, then we have
\begin{equation}\label{6l}
\left|\int_{0}^{\rho}\Lambda_6(\delta_6;w(t))dt\right|\lesssim N^{-3}\|Iw\|_{X_{1,1/2+}^{\rho}}^{6}
\end{equation}
and
\begin{equation}\label{4l}
\left|\Lambda_4(\delta_4;w(t))\right|\lesssim \|Iw\|_{H^1}^{4}. 
\end{equation}
\end{proposition}
\begin{proof}
As in \cite{C-K-S-T-T1,C-K-S-T-T3,C-K-S-T-T2}, we first perform a Littlewood-Paley decomposition of the six factors $w$, so that the $\xi_i$ are essentially the constants $N_i$, $i=1, \ldots, 6$. To recover the sum at the end we borrow a $N_s^{-\epsilon}$ from the large denominator $N_s$ and often this will not be mentioned. Also without loss of generality we can assume that the Fourier transforms in the left-side of (\ref{6l}) and (\ref{4l}) are real and nonnegative.

Let $I=\{s,a,t,b\}$ the set of indices such that $N_s\ge N_a \ge N_t \ge N_b$. We will proved first (\ref{6l}), we divide the proof into two cases.\\
1) $N_b \gtrsim N$, by definition of $m$ we have $N_t m_t \gtrsim N$ and $N_b m_b \gtrsim N$, therefore
$$N_sm^2_s\lesssim N^{-3}N_sm_s N_am_aN_tm_t N_bm_b,$$
and consequently by (\ref{int1}), H\"older inequality, (\ref{Il}) and Lemma \ref{strich}, we have
\begin{align*}
\left| \int_{0}^{\rho}\Lambda_6(\delta_6;w(t))dt\right| \lesssim &\, N^{-3}\int_{0}^{\rho}\int_{\R}\prod_{j \in I}D_xIw_j\prod_{j \notin I}w_j dx\, dt\\
\lesssim &\, N^{-3}\prod_{j \in I}\|D_xIw_j\|_{L_x^6L_{\rho}^6}\prod_{j \notin I}\|w_j\|_{L_x^6L_{\rho}^6} \\
\lesssim &\, N^{-3}\prod_{j \in I}\|Iw_j\|_{X_{1,1/2+}^{\rho}}\prod_{j \notin I}\|Iw_j\|_{X_{1-s,1/2+}^{\rho}}\\
\lesssim &\, N^{-3}\|Iw\|_{X_{1,1/2+}^{\rho}\textrm{.}}^6
\end{align*}
2) $N_b \ll N$, by (\ref{d6}) and Proposition \ref{m1}, if $N_s \ll N$, then $\Lambda_6(\delta_6)=0$, therefore we can assume $N_s \gtrsim N$, and for $\xi_1+\ldots+\xi_6=0$ this implies $N_s \sim N_a\gtrsim N$, hence
$$N_sm^2_s\lesssim N^{-1}N_sm_s N_am_a,$$
by (\ref{int1}), H\"older inequality, (\ref{Il}) and Lemmas \ref{strich} and \ref{4-2} one obtains 
\begin{align*}
\left| \int_{0}^{\rho}\Lambda_6(\delta_6;w(t))dt\right| \lesssim &\, N^{-1}\int_{0}^{\rho}\int_{\R}D_xIw_sD_xIw_a\prod_{j \notin \lbrace s, a\rbrace}w_j dx\, dt\\
\lesssim &\, N^{-1}\|(D_xIw_s) w_b\|_{L_x^4L_{\rho}^2}\|(D_xIw_a) w_p\|_{L_x^4L_{\rho}^2}\|w_t\|_{L_x^4L_{\rho}^{\infty}}\|w_q\|_{L_x^4L_{\rho}^{\infty}} \\
\lesssim &\, N^{-3}\|Iw_s\|_{X_{1,1/2+}^{\rho}}\|w_b\|_{X_{s_0,1/2+}^{\rho}}\|Iw_a\|_{X_{1,1/2+}^{\rho}}\|w_p\|_{X_{s_0,1/2+}^{\rho}} \\ 
&\|w_t\|_{X_{1/4,1/2+}^{\rho}}\|w_q\|_{X_{1/4,1/2+}^{\rho}}\\
\lesssim &\, N^{-3}\|Iw\|_{X_{1,1/2+}^{\rho}\textrm{,}}^6
\end{align*}
where $1/4<s_0<1$.

For to prove (\ref{4l}), by (\ref{delta4}) and (\ref{int1}) we have
\begin{align*}
|\Lambda_4(\delta_4;w(t))| \lesssim & \int_{\xi_1+\cdots+\xi_n}\delta_4(\xi_1,\ldots,\xi_n)\widehat{w}_1\widehat{\overline{w}}_2\widehat{w}_{3}\widehat{\overline{w}}_4\\
\lesssim & \int_{\R}|w(t)|^4dx 
\lesssim \|w(t)\|_{H^{1/4}}^4\\
\lesssim &\|Iw\|_{H^1}^{4}.
\end{align*}
Which finished the proof.
\end{proof}

%%%%%%%%%%%%%%%%%%%%%%%%%%%%%%%%%%%%%%%%%%%%%%%%%%%%%%%%%%%%%%%%
\secao{Proof of Theorem \ref{Tglobal}}
%%%%%%%%%%%%%%%%%%%%%%%%%%%%%%%%%%%%%%%%%%%%%%%%%%%%%%%%%%%%%%%%
%%%%%%%%%%%%%%%%%%%%%%%%%%%%%%%%%%%%%%%%%%%%%%%%%%%%%%%%%%%%%%%%
We will use the following results.
\begin{lemma} \label{Iu2} If $u$ is a solution of IVP (\ref{1.1}), then 
\begin{equation*}
\|Iu(t)\|_{L^2}\le \|I\varphi\|_{H^{1-s}\textrm{.}}
\end{equation*}
for $0\le s < 1$.
\end{lemma}
\begin{proof}
The lemma follows from definition of $I$, the conservation law in $L^2$ and definition of $l(\xi)$. 
%we get \begin{align}\label{Iu}
%\|Iu(t)\|_{L^2}\le & \|u(t)\|_{L^2}=\|\varphi(\xi)\|_{L^2}\le 
%\|l(\xi)\varphi(\xi)\|_{L^2}\nonumber \\
%=& \|m(\xi)\langle \xi \rangle^{1-s}\varphi(\xi)\|_{L^2\textrm{.}} 
%\end{align}
\end{proof}
\begin{lemma} \label{Iu3} If $u$ is a solution of IVP (\ref{1.1}), then 

\begin{align}\label{IT3}
|E_2(t)-E_1(t)|\le c\|I\varphi\|_{H^{1}}^{4}+cE_1(t)^{4}.%\Lambda_4(\delta_4)(t)|,
\end{align}
If $k$ is a positive integer and $u(t)$ is defined in the time interval $[0,k]$, then
\begin{align}\label{IT2}
E_2(k)=E_1(0)+ \Lambda_4(\delta_4)(0)+\sum_{j=1}^{k}\int_{j-1}^{j} \Lambda_6(\delta_6)(t) \,dt.
\end{align}
\end{lemma}
\begin{proof}
The inequality (\ref{IT3}) is obvious from (\ref{E2}), (\ref{4l}) and Lemma \ref{Iu2}.

By (\ref{d6}) we have 
\begin{align*}
E_2(k)=E_2(0)+ \sum_{j=1}^{k}\int_{j-1}^{j} \Lambda_6(\delta_6)(t) \,dt.
\end{align*}
and by (\ref{E2}) we obtain (\ref{IT2}).
\end{proof}

\subsection{Rescaling} We know that if $u(x,t)$ is a solution of (\ref{1.21}) with initial data $u(x,0)=\varphi$, then $$u_{\lambda}(x,t)=\dfrac{1}{\lambda}u(\dfrac{x}{\lambda},\dfrac{t}{\lambda^3}),$$
is also a solution of (\ref{1.21}) with initial data 
$$u_{\lambda}(x,0)=\dfrac{1}{\lambda}u(\dfrac{x}{\lambda},0)=\dfrac{1}{\lambda}\varphi(\dfrac{x}{\lambda}):=\varphi_{\lambda}.$$
Let $c_0 \in (0,1)$ a constant to be chosen later, we have 
\begin{align*}
\|I\varphi_{\lambda}\|_{H^1}\sim & \,\|\p_x I\varphi_{\lambda}\|_{L^2}+\|I\varphi_{\lambda}\|_{L^2}\\
 \lesssim &\,\dfrac{N^{1-s}}{\lambda^{1/2+s}}\|D_x^s\varphi \|_{L^2}+ \dfrac{1}{\sqrt{\lambda}}\|\varphi \|_{L^2}\\
<&\, c_0.
\end{align*}
 taking 
\begin{align}\label{lam} 
\lambda \sim N^{\textstyle\frac{2(1-s)}{1+2s}}\left(\dfrac{\|D_x^s\varphi \|_{L^2}}{c_0}\right)^{\textstyle\frac{2}{1+2s}} \quad \textrm{and} \quad N>\left(\dfrac{\|\varphi \|_{L^2}}{c_0}\right)^{\textstyle\frac{2s}{1-s}}.
\end{align}

%%%%%%%%%%%%%%%%%%%%%%%%%%%%%%%%%%%%%%%%%%%%%%%%%%%%%%
\subsection{Iteration}
%%%%%%%%%%%%%%%%%%%%%%%%%%%%%%%%%%%%%%%%%%%%%%%%%%%%%%
Without lost of generality we can assume $k_1=1$ in (\ref{E1}).
We consider our solution rescaled with initial data $$\|I\varphi \|_{H^1}=\epsilon_0<c_0<1,$$ then by Theorem \ref{tlocal} we have a solution of (\ref{1.1}) defined in the time interval $[0,1]$. For to extend the solution of local theorem in the time interval $[0,\lambda^{3}T]$ we need to prove that $\|Iu(n)\|_{H^1} \lesssim \epsilon_0$, for all $n \in \{1,2, \ldots, m_{\lambda, T}\}=W$, where $m_{\lambda, T}\sim\lambda^{3}T$ . Indeed we will prove that 
\begin{align}\label{Iter}
\|Iu(n)\|_{H^1}^2 \le 3\epsilon_0^2,\quad n \in W,
\end{align}
 but as $\|Iu(t)\|_{H^1}^2=\|Iu(t)\|_{L^2}+\|\p_xIu(t)\|_{L^2}$, by Lemma \ref{Iu2} 
 is sufficient to prove
\begin{align}\label{IT1}
\|\p_x Iu(n)\|_{L^2}^2 \le 2\epsilon_0^2, \quad n \in W.
\end{align}
We will prove (\ref{IT1}) by induction. \\
1) When $k=1$, we suppose by contradiction that $\|\p_x Iu(1)\|_{L^2}^2 >2\epsilon_0^2$, then there exist $t_0 \in (0,1)$ such that $\|\p_x Iu(t_0)\|_{L^2}^2=2\epsilon_0^2$, from (\ref{IT3}) we have
\begin{align*}
|E_2(t_0)-2\epsilon_0^2|\le 5c\epsilon_0^4,
\end{align*}
and using (\ref{E2}), (\ref{d6}), (\ref{3Iu3}) and (\ref{6l})  we obtain $$E_2(t_0)=E_1(0)+\Lambda_4(\delta_4)(0)+\int_0^{t_0}\Lambda_6(\delta_6)(0),$$ and from here
\begin{align*}%\label{IT2}
|E_2(t_0)-\epsilon_0^2|\le 5c\epsilon_0^4+\frac{1}{N^{3}}8c\epsilon_0^6,
\end{align*}
hence if $\epsilon_0^2<\frac{1}{20c}$, we have
\begin{align*}%\label{IT2}
\epsilon_0^2 \le |2\epsilon_0^2-E_2(t_0)|+|E_2(t_0)-\epsilon_0^2|\le 5c\epsilon_0^4+8c\epsilon_0^6+5c\epsilon_0^4< \epsilon_0^2,
\end{align*}
but this is a contradiction. \\
2) Now, we suppose (\ref{IT1}) for $n=1,2\ldots,k$, with $k \ge 2$ a positive integer, then we also will prove (\ref{IT1}) for $n=k+1$. In fact, in similar way as in case 1), we suppose by contradiction 
that $\|\p_x Iu(k+1)\|_{L^2}^2 >2\epsilon_0^2$, then there exist $t_0 \in (0,k+1)$ such that $\|\p_x Iu(t_0)\|_{L^2}^2=2\epsilon_0^2$. Similarly as in the case 1),
from (\ref{IT3}) we have
\begin{align}\label{C1}
|E_2(t_0)-2\epsilon_0^2|\le 5c\epsilon_0^4,
\end{align}
by (\ref{d6}) and (\ref{IT2}) we get
\begin{align*}
|E_2(t_0)-E_1(0)|\le & \left|\Lambda_4(\delta_4)(0)\right|+\left|\sum_{j=1}^{[t_0]}\int_{j-1}^{j} \Lambda_6(\delta_6)(t) \,dt\right|
+\left|\int_{[t_0]}^{t_0} \Lambda_6(\delta_6)(t) \,dt\right|,
\end{align*}
therefore by (\ref{3Iu3}) and (\ref{6l}) we easily deduce that
\begin{align}
|E_2(t_0)-\epsilon_0^2| \le & 5c\epsilon_0^4+(1+[t_0])\frac{8c}{N^{3}}\epsilon_0^6 \nonumber\\
\le &\,5c\epsilon_0^4+\lambda^{3}T\frac{8c}{N^{3}}\epsilon_0^6.
\label{C2}
\end{align}
As in the case $k=1$, by (\ref{C1}) and (\ref{C2}) we obtain a contradiction if $\lambda^{3}T \sim N^{3}$, consequently we can to iterate this process $m_{\lambda, T} \sim \lambda^{3}T$ times
if $T \sim \lambda^{-3}N^{3}$ and by (\ref{lam}) if 
$$T \sim N^{\textstyle(12s-3)/(1+2s)}.$$
Hence $u$ is globally well-posed in $H^s$ for all $s>1/4$.

\end{document}